\documentclass[11pt]{article} 
\usepackage{graphicx}
\usepackage[width=170mm,top=45mm,bottom=45mm]{geometry}
\usepackage{amssymb,amsmath} 
\usepackage{bbm}

\usepackage{mathtools} 
\usepackage{amsthm} 
\usepackage{IEEEtrantools}
\usepackage{dsfont}
\usepackage{tikz}
\usepackage{float}
\usepackage{todonotes}
\usepackage{authblk}
\usetikzlibrary{arrows.meta, positioning}
\usepackage[most]{tcolorbox}

\definecolor{hanblue}{rgb}{0.27, 0.42, 0.81}
\definecolor{mordantred19}{rgb}{0.68, 0.05, 0.0}
\definecolor{red}{rgb}{0.68, 0.05, 0.0}
\definecolor{green}{rgb}{0.0, 0.5, 0.0}

\definecolor{prooflinecolor}{gray}{0.6}

\usepackage[colorlinks, citecolor=hanblue,linkcolor=green]{hyperref}
\usepackage[capitalise]{cleveref}

\usepackage{nicefrac}
\definecolor{curveblue}{RGB}{51,153,255}

\usepackage{stmaryrd} 
\usepackage{todonotes}
\usepackage{comment}
\usepackage{color}
\usepackage{subcaption}
\usepackage{mathabx}

\theoremstyle{plain}

\newtheorem{theorem}{Theorem}[section]
\newtheorem{lemma}[theorem]{Lemma}

\newtheorem{definition}[theorem]{Definition}

\newtheorem{remark}[theorem]{Remark}

\newcommand{\ssubset}{\subset\joinrel\subset}

\newcommand{\R}{\mathbb{R}}

\newcommand{\res}{\mathop{\hbox{\vrule height 7pt width .5pt depth 0pt
\vrule height .5pt width 6pt depth 0pt}}\nolimits}

\usepackage{tikz}
\usetikzlibrary{positioning, arrows.meta, fit, calc, backgrounds}
\usepackage{caption}
\usepackage{xcolor}
\usetikzlibrary{positioning, arrows.meta}
\definecolor{boxgreen}{RGB}{204, 255, 204}      
\definecolor{bggreen}{RGB}{210, 255, 210}       
\definecolor{boxyellow}{RGB}{255, 248, 220}     
\definecolor{bordergray}{RGB}{130, 130, 130}    

\usepackage[most]{tcolorbox}

\allowdisplaybreaks

\title{A non-zero gap for the planar Mumford-Shah convex lifting and the nonexistence of crack-tip calibrations}
\author[1]{Marcello Carioni}
\author[2]{Giorgio Patrini}
\affil[1]{University of Twente, Enschede, The Netherlands}\affil[2]{Sensity, The Netherlands}

\date{}
\begin{document}

\maketitle

\begin{abstract}
We prove that the convex lifting of the Mumford–Shah functional constructed through the notion of calibration introduced in \cite{alberti2003calibration} is not sharp in dimension two. For Dirichlet boundary data prescribed by the crack-tip function, we construct a lifted competitor whose energy is strictly lower than that of the crack-tip subgraph. Since the crack-tip minimizes the original Mumford–Shah functional, this establishes a strict gap between the infima of the original and lifted problems. As a consequence we establish that the crack-tip admits no calibration.
\end{abstract}

\section{Introduction}

In this work we consider the homogeneous Mumford-Shah functional defined for functions of special bounded variation $v\in {\rm SBV}(\Omega)$ as
\begin{align*}
{\rm MS}(v) = \int_{\Omega} |\nabla v|^2 \, dx +  \mathcal{H}^{n-1}(S_v),
\end{align*}
where $\Omega \subset \R^n$ is an open, bounded set,  $S_v$ is the singular set of $v\in {\rm SBV}(\Omega)$ and $\mathcal{H}^{n-1}$ is the $n-1$ dimensional Hausdorff measure. 
Originally introduced by Mumford and Shah \cite{mumford1989optimal, mumford1985boundary} as a model for image segmentation, the full functional includes an additional fidelity term penalizing deviations from the observed image $g$:
\begin{align*}
F(v) = {\rm MS}(v) + \|v-g\|_{L^2(\Omega)}^2
\end{align*}
again considered on  $v\in {\rm SBV}(\Omega)$. The existence and regularity of minimizers of $F$, together with the characterization of local and global minimizers of ${\rm MS}$, have been the subject of extensive mathematical research over the past four decades \cite{david2005singular, de2025regularity, carioni2017existence, de2013higher, ambrosio2003higher, de2014higher}. This interest is motivated, in particular, by the central role that the characterization of minimizers of ${\rm MS}$ plays in understanding the regularity of minimizers of $F$ \cite{bonnet1996regularity}, and by its close connection to the longstanding Mumford--Shah conjecture \cite{lemenant2016selective} \footnote{At the time of writing, recent developments in this direction have appeared in \cite{deangelis2026classification} and \cite{deangelis2026solution} that, if verified, would establish the Mumford-Shah conjecture}. 
Notably, a ground-breaking result by Bonnet and David \cite{bonnet2001cracktip} has shown that in $\mathbb{R}^2$ the so-called crack-tip, defined in polar coordinates as
\begin{align}\label{eq:ca}
\hat u(\rho,\theta) = \sqrt{\frac{2\rho}{\pi}}\sin\left(\frac{\theta}{2}\right), \quad -\pi<\theta<\pi
\end{align}
is a global minimizer of  ${\rm MS}$.

With the goal of proving the minimality of candidate minimizers, in \cite{alberti2003calibration}, Alberti, Bouchitt\'e and Dal Maso introduced the theory of calibrations for the Mumford--Shah functional. Their approach builds on foundational work by Federer \cite{federer1974real} in geometric measure theory, by Harvey and Lawson \cite{harvey1982calibrated} in calibrated geometry and generalizes previous approaches related to null-lagrangians of integral functionals \cite{giaquinta2013calculus}. Calibrations are a tool to certify the minimality of a given candidate and, in case of the Mumford--Shah functional, they are divergence free vector fields $\phi$ satisfying suitable properties tailored to the candidate minimizer. In particular, their construction can provide an elegant proof of minimality with respect to either ${\rm MS}$ (with Dirichlet boundary conditions) or $F$. Following \cite{alberti2003calibration} the calibration theory for the Mumford--Shah functional and its variants has also been extensively developed in subsequent years \cite{mora2001local, dal2000local, mora2002local}.
Moreover, similar approaches have been used for various integral and free discontinuity functionals \cite{alberti2001long, pock2009convex}, free boundary problems \cite{lawlor1994paired, de2009short, carioni2020calibrations, helein1996isoperimetric, carioni2021different, pagliari2020halfspaces} and mean curvature flows \cite{fischer2022local, pluda2023minimizing}. 
It was therefore natural to seek a proof of the minimality of the crack-tip through the explicit construction of a suitable calibration. However, despite considerable efforts by several authors, including the first author of this paper, such a construction has remained elusive.

A closely related, more fundamental question is whether the crack-tip admits a calibration at all. A negative answer would explain why attempts to construct one have been unsuccessful. The main goal of this work is to show that the crack-tip does not admit a calibration of the form proposed in \cite{alberti2003calibration}.

Questions of this kind have been investigated for various classes of functionals \cite{federer1974real, bouchitte2025calibration, chambolle2001convex, carioni2016discrete}, and their answers are often closely connected to the sharpness of a suitable convex lifting of the original functional and, more generally, to the presence of a duality gap \cite{bouchitte2020convex, bouchitte2015duality}. 
For the homogeneous Mumford--Shah functional, the convex lifting has been introduced in \cite{alberti2003calibration} as a byproduct of the notion of calibration. It takes the form 
\begin{align*}
F_{{\rm lift}}(w) = \sup_{\phi \in K} \int \phi \cdot Dw 
\end{align*}
where $w \in {\rm BV}_{\rm loc}(\Omega \times \R)$ and $K$ is the set of bounded Borel vector fields $\phi: \Omega \times \R \rightarrow \R^{n+1}$ satisfying the properties: 
\begin{align*}
\left\{
\begin{array}{ll}
\phi^t(x,t) \geq \dfrac{|\phi^x(x,t)|^2}{4} \text{ for a.e.  } x\in \Omega \text{ and every } t\in \R\\
\left| \displaystyle\int_{t_1}^{t_2} \phi^x(x,s)\, ds\right|\leq 1 \text { for every } x\in \Omega \text{ and every } t_1,t_2 \in\R.
\end{array}
\right.
\end{align*}
To establish that this lifting is not sharp for prescribed boundary data, it suffices to find a Dirichlet minimizer $u$ of ${\rm MS}$ and a competitor $w$, agreeing with $\mathbf{1}_{\{u>t\}}$ (the indicator function of the subgraph of $u$) outside a given set, such that $F_{{\rm lift}}(w)<{\rm MS}(u)$, that is a non-zero gap exists. Crucially, such a non-zero gap also rules out the existence of a calibration for $u$.

In this work, we show that the convex lifting of the Mumford--Shah functional is not sharp in dimension two. We take $\hat u$ to be the crack-tip, whose minimality follows from Bonnet and David \cite{bonnet2001cracktip}. We show that for this choice a non-zero gap exists effectively ruling out the existence of calibration for the crack-tip.
Our result also identifies a limitation of numerical approaches based on this convex lifting. For boundary data exhibiting a positive gap, a minimizer of the lifted problem cannot be the subgraph lift of an admissible $SBV$ function. Therefore, solving the lifted problem does not, in general, recover a minimizer of the original functional.

We conclude this introduction by outlining the main idea of the proof. We construct a competitor $w_L\in{\rm BV}_{\rm loc}(\Omega\times\mathbb R)$ where $\Omega$ is chosen to be the annulus $\Omega = \{x \in \R^2 : e^{-4L} < |x| < 2e^{-L}\}$. Here $L$ is a positive parameter that will be chosen sufficiently big. This annular choice allows us to establish a strictly positive lifting gap in $\Omega$ and it is made to avoid the fact that calibrations for the crack-tip (if exists) need to be chosen singular close to the origin, due the the crack-tip singularity (see Remark \ref{rem:rem}). We construct $w_L$ by perturbing the lifted subgraph of the crack-tip while preserving its boundary data at the annulus. Specifically, we use the layer-cake construction
\begin{align*}
w_L(x,t) = \frac{1}{2}\int_{-1}^1  \mathbf{1}_{\{v_{L,a}>t\}}(x,t) da\quad \text{with}\quad v_{L,a}(x) = \hat u(x) + ah_L(x)
\end{align*}
where $h_L$ is a radial cut-off function vanishing for $|x|\leq e^{-3L}$ and $|x|\geq e^{-L}$ and having a plateau equal to $e^{-L}$ on a smaller annulus inside $\Omega$. 
Such perturbation is not modifying the jump set of the crack-tip.  Moreover, on the plateau, we have $\nabla v_{L,a}=\nabla\hat u$.
There, the key observation is that for any admissible field $\phi$ the lifting $\int \phi \cdot Dw_L$ cannot recover the Dirichlet energy of the crack-tip. Indeed, due to the first constraint in the definition of $K$, no-gap implies that
\begin{align*}
\phi^x(x,t)=2\nabla\hat u(x),
\qquad
\phi^t(x,t)=|\nabla\hat u(x)|^2.
\end{align*}
The  second condition in $K$ rules this out wherever $4e^{-L}|\nabla\hat u|>1$, which holds throughout a smaller annulus inside the plateau because of the crack-tip singularity: $|\nabla\hat u(x)|=(2\pi|x|)^{-1/2}$. Quantitatively, this construction produces an averaged deficit of order $e^{-2L}$, uniformly over all admissible fields.
The cutoff, that is crucially chosen logarithmically, restores the boundary data at a cost of order $e^{-2L}/L$, so choosing $L$ sufficiently large provides a competitor with strictly smaller lifted energy than the crack-tip.

\medskip

\textbf{AI statement.} The initial idea for constructing a competitor exhibiting a positive gap emerged from discussions between the authors and GPT-6 Astra in September 2026. Following several unsuccessful attempts to construct a calibration for the crack-tip, GPT-6 Astra suggested investigating the existence of a positive gap in the convex lifting and subsequently proposed the initial competitor.
The paper is built on the initial idea and it has been integrally written without the use of AI, apart from a grammar and a notational check. 
The
authors take full responsibility for the mathematical content, references
and presentation of the manuscript.

\medskip 

\textbf{Reflections and Acknowledgments \footnote{These reflections and acknowledgments concern only the first author. This is why they are written in the first person}.} Having spent years of the PhD trying to calibrate the crack-tip, this problem has obviously a deeply personal significance for me. Perhaps this explains why completing the paper and finally having an answer to this question feels bittersweet almost ten years after completing my PhD. On one hand, I am enthusiastic about the result and glad to understand why my attempts to construct a calibration could not succeed. On the other hand, knowing that the key idea came from an AI leaves me with this sense of bitterness.

Anyway I believe that some acknowledgements are in order.
First and foremost, I would like to thank my PhD advisor, Bernd Kirchheim, for introducing me to this problem and for the many hours we spent exploring strategies, working through computations, and discussing unsuccessful attempts. I am also grateful to those who worked on this problem before me and generously shared their approaches and insights. In particular, I have fond memories of discussions with Giovanni Alberti, Antonin Chambolle, Guy David, and Massimo Gobbino that also helped me get through my PhD.

\section{Preliminaries}

We first recall some basic notations and definitions we will use throughout the paper. We always denote by $\Omega$ an open, bounded set of  $\R^n$. Moreover, we denote by $\mathcal{H}^{n-1}(A)$ the $n-1$ dimensional Hausdorff measure of a given Borel set $A$. 
We will work with \emph{functions with bounded variation} defined on an open, bounded set $U$ in $\R^d$ (in this paper $d=n$ or $d = n+1$) and denote them by $\rm BV(U)$. They are defined as functions in $L^1(U)$ such that their distributional derivative is a bounded measure on $U$.  Given $v \in \rm BV(U)$ we denote by $S_v$ its singular set. Since $S_v$ is $(d-1)$-rectifiable, its measure-theoretical normal $\nu_{S_v}$ is well-defined $\mathcal{H}^{d-1}$ almost everywhere \cite{ambrosio2000functions}.
Moreover, for $\mathcal{H}^{d-1}$-a.e. $x\in S_v$ the traces of $v$, denoted by $v^+(x)$ and $v^-(x)$, are well-defined as the approximate limits of $v$ at $x \in S_v$ \cite{ambrosio2000functions}. We use the convention that $v^+(x) \geq v^-(x)$. 

It is well-known, that the derivative of a function of bounded variation $v \in {\rm BV}(U)$ can be decomposed as
\begin{align*}
Dv = \nabla v dx + (v^+ - v^-) \nu_{S_v} \mathcal{H}^{d-1}\res S_v +  D^cv
\end{align*}
where $D^c v$ is the Cantor part of the derivative \cite{ambrosio2000functions}. The set of special functions of bounded variation, denoted by ${\rm SBV(U)}$, are the functions $v \in {\rm BV}(U)$ such that $D^c v = 0$.

Given $v \in {\rm SBV}(\Omega)$ we will denote by $\textbf{1}_{\{v > t\}}$ the indicator function of its subgraph, that is 
\begin{align*}
\textbf{1}_{\{v > t\}}(x,t) = \left\{
\begin{array}{ll}
1 & \text{ if } v(x) > t\\
0 & \text{ otherwise. }
\end{array}
\right.
\end{align*}
Since $v \in {\rm SBV}(\Omega)$, it holds that $\textbf{1}_{\{v > t\}} \in {\rm SBV}(\Omega \times I)$ for every open, bounded interval $I$.

The homogeneous Mumford-Shah functional is defined for $v\in {\rm SBV}(\Omega)$ as
\begin{align*}
{\rm MS}(v) = \int_{\Omega} |\nabla v|^2 \, dx +  \mathcal{H}^{n-1}(S_v). 
\end{align*}

The study of local and global minimizers of the Mumford-Shah functional has been an area of active research since the model was defined by Mumford and Shah. Here is the precise definition that is also relevant for this paper.

\begin{definition}[Local and global minimizers]
We say that $u \in {\rm SBV}(\Omega)$ is a local minimizer of the Mumford-Shah functional if ${\rm MS}(u) < \infty$ and 
\begin{align*}
{\rm MS}(u) \leq {\rm MS}(v) 
\end{align*}
for every $v \in {\rm SBV}(\Omega)$ such that $\{u\neq v\}$ is compactly contained in $\Omega$. We also say that $u\in {\rm SBV}_{\rm loc}(\R^n)$ is a global minimizer if its restriction to any open bounded set is a local minimizer.
\end{definition}

In particular one of the most striking result concern proving the for $n=2$ the crack-tip defined in \eqref{eq:ca} is a global minimizer as stated in the next theorem due to Bonnet and David.

\begin{theorem}[Bonnet, David]
The crack-tip is a global minimizer of the Mumford-Shah functional
\end{theorem}

Due to the length of their proof, around 25 years ago, researchers have wondered if the crack-tip minimality can be proven by a calibration argument, considerably shortening it. Indeed, a calibration method was introduced by Alberti, Bouchitté and Dal Maso in 2001 \cite{alberti2003calibration} and thanks to that examples show that the construction of a calibration allowed to easily prove the local minimality of the triple junction and of harmonic functions under suitable assumptions. However, despite the effort in constructing such calibration all the attempts failed. 
We would like to summarize briefly the calibration method for the Mumford-Shah functional. 
\begin{definition}\label{def:calibration}
Given $u \in {\rm SBV}(\Omega)$, a calibration for $u$ is a $C^1$ vector field in $\phi = (\phi^x,\phi^t)  : \R^{n+1} \rightarrow \R^{n+1}$ such that ${\rm div} \phi = 0$ and with the following properties:
\begin{itemize}
\item[$a)$] $\phi^t(x,t) \geq \dfrac{|\phi^x(x,t)|^2}{4}$ for a.e. $x\in \Omega$ and every $t\in \R$,
\item[$b)$] $\left| \displaystyle\int_{t_1}^{t_2} \phi^x(x,s)\, ds\right|\leq 1$ for $\mathcal{H}^{n-1}$-a.e. $x\in \Omega$ and every $t_1,t_2 \in\R$,
\item[$c)$] $\phi^x(x,u(x)) = 2\nabla u(x)$ and $\phi^t(x,u(x)) = |\nabla u(x)|^2$ for a.e. $x\in \Omega$,
\item[$d)$] $ \displaystyle\int_{u^-(x)}^{u^+(x)} \phi^x(x,s)\, ds= \nu_{S_u}(x)$ for $\mathcal{H}^{n-1}$-a.e. $x\in S_u$.
\end{itemize}
\end{definition}
As shown in \cite{alberti2003calibration} the existence of a calibration of a given $u \in {\rm SBV}(\Omega)$ automatically implies its local minimality in $\Omega$. Moreover, it is worth to remark that the $C^1$ regularity of $\phi$ can be relaxed by essentially asking enough regularity to ensure the validity of the divergence theorem. For example, in \cite{alberti2003calibration}, the class of approximately regular vector fields has been considered allowing for discontinuous calibrations.

\subsection{The convex lifting for the Mumford-Shah functional and its gap}

Due to the presence of the $\mathcal{H}^{n-1}(S_v)$ the Mumford-Shah functional is naturally non-convex. However, the set of calibrations can be used to define a convex lifting. More precisely, for $w \in {\rm BV}_{\rm loc}(\Omega \times \R)$ one can define the convex functional
\begin{align*}
F_{{\rm lift}}(w) = \sup_{\phi \in K} \int \phi \cdot Dw 
\end{align*}
where $K$ is the set of bounded Borel vector fields $\phi: \Omega \times \R \rightarrow \R^{n+1}$ satisfying the properties: 
\begin{align*}
\left\{
\begin{array}{ll}
\phi^t(x,t) \geq \dfrac{|\phi^x(x,t)|^2}{4} \text{ for a.e.  } x\in \Omega \text{ and every } t\in \R\\
\left| \displaystyle\int_{t_1}^{t_2} \phi^x(x,s)\, ds\right|\leq 1 \text { for every } x\in \Omega \text{ and every } t_1,t_2 \in\R.
\end{array}
\right.
\end{align*}
It is easy to see that for every $v\in {\rm SBV}(\Omega)$ it holds that  $F_{{\rm lift}}(\textbf{1}_{\{v > t\}}) = {\rm MS}(v)$ (Remark 3.3 \cite{alberti2003calibration}). 
However, a much deeper question is establishing if the convex lifting is sharp, meaning if it has the same minimum value as the Mumford-Shah functional under suitable boundary conditions. In more precise words, given $u \in {\rm SBV}(\Omega)$ a local minimizer of ${\rm MS}$ one asks if  
\begin{align*}
\inf_{w \in \mathcal{G}_{u}} F_{\rm lift}(w) = \rm MS(u).  
\end{align*}
where $\mathcal{G}_{u} = \{v \in {\rm BV}_{\rm loc}(\Omega \times \R) : \{v \neq \textbf{1}_{\{u >t\}}\}  \ssubset \Omega \times \R\}$.
In dimension one ($n=1$), this was proven to be true by Chambolle in \cite{chambolle2001convex}. The goal of this paper is to show that this is false in dimension two. 

The lifting $F_{{\rm lift}}$ is connected to the notion of calibration by the following simple observation. Given $u \in {\rm SBV}(\Omega)$ a local minimizer of the Mumford-Shah functional, if there exists a calibration for $u$ according to Definition \ref{def:calibration}, then 
\begin{align*}
\inf_{w \in \mathcal{G}_{u}} F_{\rm lift}(w) = {\rm MS}(u).  
\end{align*}
This can be easily seen by observing that if $\bar \phi$ is a calibration for $u$, then for all $w \in {\rm BV}_{\rm loc}(\Omega \times \R)$ agreeing with $\textbf{1}_{\{u>t\}}$ outside a compact set contained in $\Omega \times \R$, it holds 
\begin{align*}
F_{\rm lift}(w) \geq \int_{\Omega \times \R} \bar \phi \cdot Dw = \int_{\Omega \times \R} \bar \phi \cdot D\textbf{1}_{\{u > t\}} = {\rm MS}(u).
\end{align*}
In this paper we show that this is not the case when $\hat u$ is the crack-tip.

\section{Non-zero gap for the crack-tip}
In this section, let $L>0$ be a parameter that will be chosen big enough and we consider the annulus
\begin{align*}
\Omega =\{x \in \R^2 : e^{-4L} < |x| < 2e^{-L}\}.
\end{align*}
Note that for notational convenience we do not make the dependence on $L$ in $\Omega$ explicit.
Moreover, we set $I = (-2,2)$ and we let $\hat u \in {\rm SBV}(\Omega)$ be the crack-tip. Our goal is to show that for big enough $L$ it holds that 
\begin{align*}
\inf_{w \in \mathcal{G}_{\hat u}} F_{\rm lift}(w) < {\rm MS}(\hat u)
\end{align*}
where 
\begin{align*}
 \mathcal{G}_{\hat u} = \{w \in {\rm BV}(\Omega \times I) : \{w \neq \textbf{1}_{\{\hat u >t\}}\}  \ssubset \Omega \times I\}.
\end{align*}
In particular, we aim to construct a competitor $w \in \mathcal{G}_{\hat u}$ such that $F_{\rm lift}(w) < {\rm MS}(\hat u)$. By proving this statement we will show the non-existence of calibration for the crack-tip in $\Omega$ and also in any set $\Omega' \supset \Omega$  by a simple restriction argument. Also note that the set $I$ is useful to work with functions defined on compact cylinders. It is enough to choose any open interval containing the image of the crack-tip in $\Omega$. From now on we will work with functions defined in $\Omega \times I$, replacing also ${\rm BV}_{\rm loc}$ with ${\rm BV}$. 

\begin{remark}\label{rem:rem}
Note that in this paper we always work with the annulus $\Omega$ as a domain. This choice is done to have a non-trivial result. Indeed, due to the condition $\phi^x = 2\nabla \hat u$ in Definition \ref{def:calibration} any calibrating vector field would have the singularity of the crack-tip at the origin. Therefore there is no hope to find $C^1$ calibrations according to Definition \ref{def:calibration}, or even approximately regular   calibrations according to \cite{alberti2003calibration}.
\end{remark}
Our goal is to construct a family of functions $w_L : \Omega \times I \rightarrow \R$ and to show that for suitable choice of $L>0$ such $w_L$ belongs to $\mathcal{G}_{\hat u}$ and $F_{\rm lift}(w_L) < {\rm MS}(\hat u)$. We first define
\begin{align*}
q(\tau) = \left\{\begin{array}{ll}
1 & |\tau| \leq \frac{1}{2}\\
2(1-|\tau|) & \frac{1}{2}<|\tau|<1\\
 0 & |\tau|\geq 1
\end{array}
\right.
\end{align*}
and the cut-off function $h_L : \Omega \rightarrow \R$
\begin{align*}
h_L(x) = e^{-L} q\left(\frac{\log(|x|) + 2L}{L}\right).
\end{align*}
The family of competitors $w_L: \Omega \times I \rightarrow \R$ is defined as:
\begin{align*}
v_{L,a}(x) = \hat u(x) + ah_L(x) \ \  \ a \in (-1,1), \qquad w_{L}(x,t) = \frac{1}{2}\int_{-1}^1  \mathbf{1}_{\{v_{L,a}>t\}}(x,t) da.
\end{align*}

It is easy to see that for $L \geq 2$ the function $w_L$ is indeed an admissible competitor as next lemma shows.
\begin{lemma}
The function $w_L$ belongs to ${\rm BV}(\Omega \times I)$ and, for $L \geq 2$, it holds that
\begin{itemize}
\item[i)]$w_L(x,t) = \mathbf{1}_{\{\hat u>t\}}(x,t)$ if $x \in \Omega$, $|x| \leq e^{-3L}$ or $|x| \geq e^{-L}$, and $t\in I$,
\item[ii)] $w_L(x,t) = 1$ for all $(x,t) \in \Omega \times (-2,-1)$,
\item[iii)] $w_L(x,t) = 0$ for all $(x,t) \in \Omega \times (1,2)$.
\end{itemize}
In particular $w_L\in \mathcal{G}_{\hat u}$.
\end{lemma}

\begin{proof}
First, let us prove that $w_L \in {\rm BV}(\Omega \times I)$. Note that $w_L$ is bounded in $\Omega \times I$ and thus $L^1(\Omega \times I)$. Moreover, since the total variation of its distributional derivative can be computed as 
\begin{align*}
|Dw_L|(\Omega\times I) & = \frac{1}{2} \int_{-1}^1 |D\textbf{1}_{v_{L,a}}|(\Omega\times I) da = \\
& = \frac{1}{2} \int_{-1}^1 \left(\int_\Omega \sqrt{1 + |\nabla \hat u(x) + a \nabla h_L(x)|^2} dx + \int_{S_{\hat u}}|{\hat u}^+(x) - {\hat u}^-(x)| d\mathcal{H}^1(x)\right) da \\
& < \infty,
\end{align*}
since $h_L$ is Lipschitz. 
Since $q(\tau)=0$ for $|\tau|\geq1$, we have
$h_L(x)=0$ whenever $|x|\leq e^{-3L}$ or
$|x|\geq e^{-L}$. This proves (i). Note also that for $x \in \Omega$, $a\in (-1,1)$ and $L \geq 2$ it holds  
\begin{align*}
v_{L,a}(x) = \hat u(x) + ah_L(x)  \leq \sqrt{\frac{2}{\pi}} + e^{-L} \leq  \sqrt{\frac{2}{\pi}} + e^{-2} < 1
\end{align*}
proving $iii)$. $ii)$ can be proven similarly, implying that $w_L\in \mathcal{G}_{\hat u}$.
\end{proof}
We also compute the Dirichlet energy of $h_L$ in $\Omega$.
\begin{lemma}\label{lem:dirichlet2}
It holds that 
\begin{align*}
\int_{\Omega} |\nabla h_L(x)|^2\, dx = e^{-2L} \frac{8 \pi}{L}.
\end{align*}
\end{lemma}
\begin{proof}
We compute
\begin{align*}
\int_{\Omega} |\nabla h_L(x)|^2\, dx & = \frac{e^{-2L}}{L^2} \int_{\Omega} q'\left(\frac{\log(|x|) + 2L}{L}\right)^2\frac{1}{|x|^2}\, dx  = \frac{2\pi e^{-2L}}{L^2} \int_{e^{-4L}}^{2e^{-L}} q'\left(\frac{\log(r) + 2L}{L}\right)^2 \frac{1}{r} dr .
\end{align*}
Set now $\tau = \frac{\log(r) + 2L}{L}$. This gives
\begin{align*}
= \frac{2\pi e^{-2L}}{L} \int_{-2}^{\frac{\log 2 + L}{L}} q'\left(\tau\right)^2 d\tau =   \frac{2\pi e^{-2L}}{L} \left(\int_{-1}^{-1/2} 4 d\tau + \int_{1/2}^1 4 d\tau\right) =   e^{-2L}\frac{8\pi}{L}
\end{align*}
as we wanted to show.
\end{proof}

\begin{theorem}\label{thm:main}
There exists a sufficiently large $L > 2 \log 4$ such that 
\begin{align*}
F_{\rm lift}(w_L)  < {\rm MS}(\hat u).
\end{align*}
\end{theorem}

\begin{proof}
Let  $\phi = (\phi^x, \phi^t) \in K$ and $v \in SBV(\Omega)$ with $|v(x)| \leq 1$ almost everywhere. Then denoting by $Q = \Omega \times I$, we can compute
\begin{align*}
{\rm MS}(v) & - \int_Q \phi \cdot D\textbf{1}_{\{v(x)>t\}} =  \int_\Omega |\nabla v(x)|^2 - \phi^x(x,v(x)) \cdot \nabla v(x)  + \phi^t(x,v(x)) \, dx \\
& \quad + \left(\int_{S_v} 1 - \nu_{S_v}(x) \cdot \int_{v^-(x)}^{v^+(x)} \phi^x(x,t)\, dt\right) d\mathcal{H}^1(x)\\
& = \int_\Omega \left| \nabla v(x) - \frac{1}{2} \phi^x(x,v(x))\right|^2 + \phi^t(x,v(x)) - \frac{1}{4} \left|\phi^x(x,v(x))\right|^2\, dx +\\
& \quad  + \left(\int_{S_v} 1 - \nu_{S_v}(x) \cdot \int_{v^-(x)}^{v^+(x)} \phi^x(x,t)\, dt\right) d\mathcal{H}^1(x).
\end{align*}
Applying the previous computation with $v= v_{L,a}$ it holds that 
\begin{align}\label{eq:one}
& \int_Q  \phi \cdot Dw_L = \frac{1}{2} \int_Q\int_{-1}^1 \phi \cdot D \textbf{1}_{\{v_{L,a}(x) > t\}}\, da \nonumber\\
& = \frac{1}{2} \int_{-1}^1 \int_Q  \phi \cdot D \textbf{1}_{\{v_{L,a}(x) > t\}}\, da - \frac{1}{2} \int_{-1}^1 {\rm MS}(v_{L,a}) da + \frac{1}{2} \int_{-1}^1 {\rm MS}(v_{L,a}) da \nonumber\\
&  = \frac{1}{2} \int_{-1}^1 {\rm MS}(v_{L,a}) da - \frac{1}{2}\int_{-1}^1\int_\Omega \left| \nabla v_{L,a}(x) - \frac{1}{2} \phi^x(x,v_{L,a}(x))\right|^2 + \phi^t(x,v_{L,a}(x)) - \frac{1}{4} \left|\phi^x(x,v_{L,a}(x))\right|^2\, dx da \nonumber \\
& \qquad - \frac{1}{2}\int_{-1}^1 \left(\int_{J_{v_{L,a}}} 1 - \nu_{v_{S_{L,a}}}(x) \cdot \int_{v_{L,a}^-(x)}^{v_{L,a}^+(x)} \phi^x(x,t)\, dt\right) d\mathcal{H}^1(x)da \nonumber \\
& \leq \frac{1}{2} \int_{-1}^1 {\rm MS}(v_{L,a}) da - \frac{1}{2}\int_{-1}^1\int_\Omega \left| \nabla v_{L,a}(x) - \frac{1}{2} \phi^x(x,v_{L,a}(x))\right|^2 dx da
\end{align}
where in the last inequality we used that $\phi \in K$.
Now note that 
\begin{align}\label{eq:two}
\frac{1}{2} & \int_{-1}^1 {\rm MS}(v_{L,a}) da = \frac{1}{2} \int_{-1}^1 \int_{\Omega} |\nabla \hat u(x) + a \nabla h_L(x)|^2\,  dx\, da  + \frac{1}{2}\int_{-1}^1 \mathcal{H}^1(S_u) \nonumber \\
&  = {\rm MS}(u) + \frac{1}{2} \int_{-1}^1 \int_{\Omega} a^2 |\nabla h_L(x)|^2   dx \, da +   \int_{-1}^1 \int_{\Omega} a \nabla \hat u(x) \cdot \nabla h_L(x)   dx \, da  \nonumber \\
& =  {\rm MS}(u) + \frac{1}{2} \int_{-1}^1 \int_{\Omega} a^2 |\nabla h_L(x)|^2   dx \, da \nonumber \\
& = {\rm MS}(u) + \frac{1}{3} \int_{\Omega} |\nabla h_L(x)|^2   dx.
\end{align}
Moreover, using Jensen's inequality and the fact that $\phi \in K$ one has that 
\begin{align}\label{eq:three}
 \frac{1}{2}\int_{-1}^1 \int_\Omega & \left| \nabla v_{L,a}(x) - \frac{1}{2} \phi^x(x,v_{L,a}(x)) \right|^2 \, dx \, da \nonumber \\
 & = \frac{1}{2}\int_\Omega \int_{-1}^1 \left| \nabla \hat u(x) + a \nabla h_L(x) - \frac{1}{2} \phi^x(x, u(x) + a h_L(x))\right|^2\, da dx \nonumber \\
& \geq \int_\Omega  \left| \frac{1}{2}\int_{-1}^1 \nabla \hat u(x) + a \nabla h_L(x) - \frac{1}{2} \phi^x(x, \hat u(x) + a h_L(x)) da \right|^2 dx   \nonumber \\
& = \int_\Omega  \left| \frac{1}{2} \int_{-1}^1 \nabla \hat u(x)  - \frac{1}{2} \phi^x(x, \hat u(x) + a h_L(x)) da \right|^2 dx \nonumber\\
& = \int_\Omega  \left| \nabla \hat u(x)  -  \int_{-1}^1 \frac{1}{4} \phi^x(x, \hat u(x) + a h_L(x)) da \right|^2 dx \nonumber \\
& \geq \int_{\Omega \cap \{h_L >0\}}  \left| \nabla \hat u(x)  -  \frac{1}{4h_L(x)} \int_{\hat u(x) - h_L(x)}^{\hat u(x) + h_L(x)}  \phi^x(x, z) dz \right|^2 dx  \nonumber \\
& \geq \int_{\Omega \cap \{h_L >0\}} \left( |\nabla \hat u(x)|  -  \frac{1}{4h_L(x)} \right)^2_+ dx 
\end{align}
for all $\phi \in K$. Therefore defining the annulus
\begin{align*}
\mathcal{A}_L = \left\{x: \frac{e^{-2L}}{4} \leq |x| \leq e^{-2L}\right\} 
\end{align*} 
we note that $h_L(x) = e^{-L}$ in $\mathcal{A}_L$. Indeed, since $L > 2 \log 4$, we get that $0 \geq \frac{\log(|x|) + 2L}{L} \geq \frac{-\log 4}{L} > -\frac{1}{2}$. Therefore, using Lemma \ref{lem:dirichlet2} 
\begin{align*}
\int_Q  \phi \cdot Dw_L  & \leq {\rm MS}(\hat u) + \frac{1}{3} \int_{\Omega} |\nabla h_L(x)|^2   dx -  \int_{\Omega \cap \{h_L >0\}}   \left( |\nabla \hat u(x)|  -  \frac{1}{4h_L(x)} \right)^2_+ dx  \\
& \leq {\rm MS}(\hat u)+  e^{-2L} \frac{8 \pi}{3L} - \int_{\mathcal{A}_L}  \left( \frac{1}{\sqrt{2\pi |x|}} -  \frac{1}{4e^{-L}} \right)^2_+ dx\\
& =  {\rm MS}(\hat u)+  e^{-2L} \frac{8 \pi}{3L} - \int_{\mathcal{A}_L}  \left( \frac{1}{\sqrt{2\pi |x|}} -  \frac{1}{4e^{-L}} \right)^2 dx\\
& = {\rm MS}(\hat u)+  e^{-2L} \frac{8 \pi}{3L} - 2\pi e^{-2L} \int_{1/4}^1 r\left( \frac{1}{\sqrt{2\pi r}} -  \frac{1}{4} \right)^2 dr\\
& = {\rm MS}(\hat u)+  e^{-2L} \left(\frac{8 \pi}{3L} - 2\pi  \int_{1/4}^1 r\left( \frac{1}{\sqrt{2\pi r}} -  \frac{1}{4} \right)^2 dr \right).
\end{align*}
Therefore, choosing $L$ big enough and since the right hand side does not depend on $\phi$ we immediately obtain that $\sup_{\phi \in K} \int_Q  \phi \cdot Dw_L < {\rm MS}(\hat u)$ as we wanted to prove.

\begin{remark}
Our construction provides an explicit lower bound on the
relaxation gap, although it does not determine its exact value.
By considering ${\rm MS}$ and $F_{\mathrm{lift}}$ defined on $\Omega_L =\{x \in \R^2 : e^{-4L} < |x| < 2e^{-L}\}$ we define
\[
\Delta_L
:=
\mathrm{MS}(\hat u)
-
\inf_{w\in\mathcal G_{\hat u}}F_{\mathrm{lift}}(w).
\]
Set
\begin{align*}
\eta
&:=
2\pi\int_{1/4}^{1}
r\left(\frac{1}{\sqrt{2\pi r}}-\frac{1}{4}\right)^2 dr =
\frac{3}{4}-\frac{7\sqrt{2\pi}}{24}+\frac{15\pi}{256}
\approx 0.203
\end{align*}
The estimate in the proof of Theorem \ref{thm:main} thus gives
\begin{align*}
\Delta_L
\geq 
e^{-2L}\left(\eta-\frac{8\pi}{3L}\right).
\end{align*}
This lower bound is strictly positive whenever $L> \frac{8\pi}{3\eta} \approx 41.273$. So choosing $L = 42$ we get $\Delta_{42} > e^{-90}$. This is obviously a crude estimate; however, we believe that the real gap is also small. This could explain why a positive gap has never been convincingly detected using numerical methods.
\end{remark}

\end{proof}

\footnotesize

\bibliographystyle{abbrv}

\bibliography{biblio}

\end{document}